\documentclass[11pt,reqno]{amsart}

\usepackage[top=1in, bottom=1in, left=1.3in, right=1.3in]{geometry}
\usepackage{changes}

\usepackage{amsmath}
\usepackage{amssymb}
\usepackage{amsthm}
\usepackage{verbatim}
\usepackage{graphicx}
\usepackage{graphics}
\usepackage{color}
\usepackage{enumerate}
\usepackage{mathrsfs}
\usepackage{booktabs}
\usepackage[toc,page]{appendix}
\usepackage{enumerate}
\usepackage{mathrsfs}
\usepackage{fancyhdr}
\usepackage{amsopn, amstext, amscd,pifont}
\usepackage{amssymb,bm, cite}
\usepackage[all]{xy}
\usepackage{color}
\usepackage{graphicx,hyperref,float}
\usepackage{float}
\usepackage{listings}
\usepackage{setspace}
\usepackage{url}
\usepackage{paralist}
\usepackage{subcaption}
\usepackage{stmaryrd}
\usepackage{tikz-cd}

\usepackage{graphicx,pgf}
\usepackage{float}
\usepackage{enumitem}

\usepackage{fancyhdr}

\numberwithin{equation}{section} 

\newtheorem{introtheorem}{Theorem}

\theoremstyle{plain}

\newtheorem{theorem}{Theorem}[section]
\newtheorem{proposition}[theorem]{Proposition}
\newtheorem{lemma}[theorem]{Lemma}

\newtheorem*{conjecture*}{Conjecture}

\theoremstyle{definition}
\newtheorem{definition}[theorem]{Definition}

\graphicspath{{images/}}

\newcommand{\Crit}{\operatorname{Crit}}
\newcommand{\diam}{\operatorname{diam}}
\newcommand{\Rep}{\operatorname{Rep}}
\newcommand{\Dom}{\operatorname{Dom}}
\newcommand{\dist}{\operatorname{Dist}}

  {\end{list}}

\title[Repelling Density for Algebraic Rational Maps]{Repelling Density for Algebraic Rational Maps}
\author{Shilei Fan}
\address{School of Mathematics and Statistics, and Key Laboratory of Nonlinear Analysis \& Applications (Ministry of Education), Central China Normal University, Wuhan 430079, P. R. China}  \email{slfan@mail.ccnu.edu.cn}

\author{Hongming Nie}
\address{Los Angeles, CA 90089, USA}
\email{hongming.i.nie@gmail.com}

\author{Yuefei Wang}
\address{School of Mathematical Sciences, Shenzhen University, Shenzhen 518060, Guangdong, China \& Academy of Mathematics and System Sciences, CAS, Beijing 100190, China}
\email{wangyf@math.ac.cn}

\thanks{S. L. FAN was partially supported by NSFC (grants No.  12331004).  Y. F. WANG was partially supported by NSFC (grants No. 12231013)}

\begin{document}

\begin{abstract}
We prove that, for every rational map of degree at least $2$ with  algebraic  coefficients over $\mathbb{Q}_p$, the  classical Julia set  in $\mathbb{P}^1(\mathbb{C}_p)$  coincides with the closure of its repelling periodic points. This establishes the Repelling Density Conjecture for rational maps with algebraic  coefficients over $\mathbb{Q}_p$.
\end{abstract}
\maketitle

\section{Introduction}\label{sec:intro}

The iteration theory of rational maps on the Riemann sphere originated in the work of Fatou and Julia in the early twentieth century. Its basic dichotomy separates the Fatou set, where the iterates form a normal family, from its complement, the Julia set. A fundamental theorem asserts that the Julia set of every rational map of degree at least two is the closure of its repelling periodic points. Fatou and Julia established this theorem by different methods; both proofs are presented in \cite[\S 14]{Milnor06}. Fatou's argument combines non-normality with the finiteness of nonrepelling cycles. Julia's argument starts from a fixed point in the Julia set and constructs a homoclinic backward orbit passing through a prescribed open set that meets the Julia set, and then employs a suitable return construction to yield a repelling periodic point in that open set.

The corresponding density problem in non-archimedean dynamics is more subtle. For the non-archimedean dynamics, we refer the reader to \cite{Benedetto19}. In this paper, we focus on $p$-adic rational maps.  Let $\mathbb C_p$ be the completion of an algebraic closure of $\mathbb Q_p$, equipped with  the $p$-adic absolute value $|\cdot|$. For a rational map $f\in\mathbb C_p(z)$ of degree at least two, its classical Fatou set $F(f)$ consists of the points of $\mathbb P^1(\mathbb C_p)$ having a neighborhood on which the iterates of $f$ are equicontinuous with respect to the spherical metric. Its complement is the classical Julia set $J(f)$. A periodic point in $\mathbb P^1(\mathbb C_p)$ is called \emph{repelling} if its multiplier has absolute value greater than one. We denote the set of repelling periodic points by $\operatorname{Rep}(f)$. 

The following conjecture is the $p$-adic version of the general non-archimedean Repelling Density Conjecture. The conjecture goes back to Hsia \cite{Hsia00} and  is stated in \cite[Conjecture 17.1]{BenedettoS19}. 

\begin{conjecture*}
Let $f\in\mathbb C_p(z)$ be a rational map of degree at least two. Then
$$J(f)=\overline{\operatorname{Rep}(f)}.$$
\end{conjecture*}

Adapting Fatou's method, Hsia \cite{Hsia00} proved that the classical Julia set is contained in the closure of all periodic points. Hsia's result does not imply repelling density, since a non-archimedean rational map may have infinitely many non-repelling cycles.  B\'ezivin \cite{Bezivin01} adapted Julia's method and proved that the existence of a single repelling periodic point implies the conjectured density. Thus the conjecture reduces to showing that every nonempty classical Julia set contains a repelling periodic point.

The conjecture has been established in several cases. Briend and Hsia \cite{Briend12} proved it for cubic polynomials. Doyle, Jacobs and Rumely \cite{Doyle16} proved it for quadratic rational maps. Favre and Rivera-Letelier \cite{Favre26} established it when $p>\deg f$. More recently, Truong \cite{Truong26} proved it for cubic rational maps, using an approach related to an idea of T. Silverman \cite{Silverman20}.

In this paper, we confirm the conjecture for rational maps with algebraic coefficients. 

\begin{introtheorem}\label{main}
Let $f\in\overline{\mathbb{Q}}_p(z)$ be a rational map of degree at least two. Then
$$J(f)=\overline{\operatorname{Rep}(f)}.$$
\end{introtheorem}

Theorem \ref{main} follows from the following criterion, which imposes no algebraicity assumption on the coefficients.

\begin{introtheorem}\label{mainB} 
Let $f\in\mathbb C_p(z)$ be a rational map of degree at least two. If $J(f)$ contains a nonempty compact forward-invariant subset, then
$$J(f)=\overline{\operatorname{Rep}(f)}.$$
\end{introtheorem}

 In view of B\'ezivin's theorem  \cite{Bezivin01}, the essential assertion of Theorem \ref{mainB} is that the existence of such a subset forces the existence of a repelling periodic point. 

To establish Theorem \ref{mainB}, we construct suitable puzzles and consider the return maps, which can be recorded by a natural collection of multigraphs. We then relate the closed walks in these multigraphs to the existence of a repelling periodic cycle. To deduce Theorem \ref{main} from Theorem \ref{mainB}, we show that  any rational map $f\in\overline{\mathbb{Q}}_p(z)$ of degree at least $2$ with nonempty $J(f)$ has an algebraic point over $\mathbb Q_p$ in $J(f)$, whose forward orbit has a compact and forward-invariant closure in $J(f)$. 

\medskip
\textbf{AI usage:} We discuss the conjecture with ChatGPT for rational maps in $\mathbb{C}_p(z)$ that has a compact forward-invariant subset of the Julia set,  building on the framework  established in previous works \cite{Fan21,Fan24}  and \cite{Hsia25}. The authors  developed  the main idea and the ChatGPT helped clarify some details. In particular, the derivative identities in \S\ref{sec:lemmas} are suggested by ChatGPT. All mathematical statements and proofs are checked and finalized by the authors. 

\section{Auxiliary Lemmas}\label{sec:lemmas}
In this section, we prove several lemmas about the derivatives of rational maps $f\in\mathbb{C}_p(z)$. For a disk $D\subseteq \mathbb{C}_p$ containing no poles of $f$,  denote by $||f||_D$ the Gauss norm of $f$ on $D$, that is $||f||_D:=\sup\{|f(x)|: x\in D\}$.  Moreover, we write $\overline{D}(a,r):=\{z\in\mathbb{C}_p: |z-a|\le r\}$ (resp. $D^-(a,r):=\{z\in\mathbb{C}_p: |z-a|<r\}$) for the closed (resp. open) disk centered at $a\in \mathbb{C}_p$ with radius $r>0$, where $r=\diam(\overline{D}(a,r))=\diam(D^-(a,r))$ is the diameter. 

We begin with the following lemma.
\begin{lemma}\label{lem:basic}
Let $D\subseteq\mathbb{C}_p$ be a disk, and let $f\in\mathbb{C}_p(z)$ be nonconstant and have no pole on $D$. Then $f(D)$ is a disk, and for all $x,y\in D$,
$$|f(x)-f(y)|\le\frac{\diam(f(D))}{\diam(D)}\cdot|x-y|.$$
In particular,  for any $z\in D$, 
$$|f'(z)|\le\frac{\diam(f(D))}{\diam(D)},$$
and  the equality holds if the Weierstrass degree of $f$ on $D$ is $1$. 
Moreover, if $f$ has no critical point in $D$, $|f'(z)|$ is constant on $D$. 
\end{lemma}
\begin{proof}
The assertion, except the ``Moreover" part,  follows from \cite[Propositions 3.20 and 3.25]{Benedetto19}. For the  ``Moreover" part, note that $f'(D)$ is a disk with $0\not\in f'(D)$. Then $\diam(f'(D))\le|a|$ for any $a\in f'(D)$. Thus, for any $b\in f'(D)$, 
$$|a-b|\le \diam(f'(D))\le\min\{|a|, |b|\},$$ 
 and hence $|a|=|b|$. 
\end{proof}

We now establish two identities concerning the derivatives and distances $|f(c)-c|$ for critical points $c$ of $f$. The first identity (Lemma \ref{lem:fix}) compares fixed-point multipliers with $|f(c)-c|$. The second (Lemma \ref{lem:any}) applies a similar argument to an arbitrary point and its preimages.

\begin{lemma}\label{lem:fix}
Let $a\in\mathbb C_p$ and $r,R\in|\mathbb C_p^\times|$ satisfy $r<R$. 
Let $f\in\mathbb C_p(z)$ be a rational map with no poles in the closed disk $\overline{D}:=\overline D(a,r)$ and satisfying $f(\overline D)=\overline D(a,R)$.
Then
$$\prod_{\alpha\in \overline{D}: f(\alpha)=\alpha}|f'(\alpha)|=(||f'||_{\overline D})^{|\{\alpha\in \overline{D}: f(\alpha)=\alpha\}|}\prod_{c\in \overline{D}: f'(c)=0}\frac{|f(c)-c|}{R},$$
where both fixed points and critical points are  counted with multiplicity, and empty products take value $1$. 
\end{lemma}

\begin{proof}
Set $g(z):=f(z)-z.$
Since $f$ has no poles on $\overline D$, both $g$ and $f'$ are analytic on $\overline D$. For $x\in \overline D(0,r)$, write $$f(a+x)=f(a)+\sum_{j\ge1}b_jx^j.$$ 
Then 
$$g(a+x)=f(a)-a+(b_1-1)x+\sum_{j\ge2}b_jx^j.$$

Since $f(\overline{D})=\overline D(a,R)$, we have $|f(a)-a|\le R$  and $\max_{j\ge1}|b_j|r^j=R$. Denote 
$$q:=\max\{j\ge1:|b_j|r^j=R\}.$$
Since  $r<R$, we have 
$$|b_q|r^q=R\ge\max\{|f(a)-a|, |(b_1-1)|r\}.$$
 Hence  $||g||_{\overline D}=R$ and the  Weierstrass degree of $g$ on $\overline D$ is $q$. 
It follows that $g$ has exactly $q$ zeros, say $\alpha_1,\dots,\alpha_q$, in $\overline D$, counted with multiplicity.

By Weierstrass preparation, we may write
$$g(z)=U(z)\prod_{i=1}^{q}(z-\alpha_i), \qquad f'(z)=V(z)\prod_{j=1}^{m}(z-c_j),$$
where $U$ and $V$ are analytic units on $\overline D$, and $c_j\in \overline D$ with $f'(c_j)=0$. If there is no such $c_j$, we simply write $f'(z)=V(z)$ and take $m=0$.  
Since such units have constant absolute value on $\overline D$, we have that, for any $z\in \overline D$, 
$$|U(z)|=\frac{R}{r^q}, \qquad |V(z)|=\frac{||f'||_{\overline D}}{r^m}.$$

Evaluating $f'$ at the fixed points yields
$$\prod_{i=1}^{q}|f'(\alpha_i)|=\frac{(||f'||_{\overline D})^q}{r^{mq}}\prod_{i=1}^{q}\prod_{j=1}^{m}|\alpha_i-c_j|, $$
and evaluating $g(z)$ at the critical points gives
$$\prod_{j=1}^{m}|f(c_j)-c_j|=\prod_{j=1}^{m}|g(c_j)|=\frac{R^m}{r^{qm}}\prod_{j=1}^{m}\prod_{i=1}^{q}|c_j-\alpha_i|,$$
where if $m=0$, products take value $1$.
It follows that 
$$\prod_{i=1}^{q}|f'(\alpha_i)|=(||f'||_{\overline D})^q\prod_{j=1}^{m}\frac{|f(c_j)-c_j|}{R}.$$
 The conclusion holds.
\end{proof}

\begin{lemma}\label{lem:any}
Let $a,b \in\mathbb C_p$ and $r,R\in|\mathbb C_p^\times|$.  
Let $f\in\mathbb C_p(z)$ be a rational map  that maps $\overline{D}:=\overline D(a,r)$ onto $\overline D(b,R)$ with degree $d\ge 1$.
Then for any $z_0\in \overline D(b,R)$, 
$$\prod_{\alpha\in \overline{D}: f(\alpha)=z_0}|f'(\alpha)|=(||f'||_{\overline D})^d\prod_{c\in \overline{D}: f'(c)=0}\frac{|f(c)-z_0|}{R},$$
where preimages and critical points are  counted with multiplicity, and empty
products take value $1$. 
\end{lemma}

\begin{proof}
By Weierstrass preparation, we may write
$$f(z)-z_0=U(z)\prod_{i=1}^{d}(z-\alpha_i), \qquad f'(z)=V(z)\prod_{j=1}^{m}(z-c_j),$$
where $U$ and $V$ are analytic units on $\overline D$, $\alpha_i\in \overline D$ with $f(\alpha_i)=z_0$, and $c_j\in \overline D$ with $f'(c_j)=0$. If there is no such $c_j$, we simply write $f'(z)=V(z)$ and take $m=0$.  
Since such units have constant absolute value on $\overline D$, we have that, for any $z\in \overline D$, 
$$|U(z)|=\frac{R}{r^d}, \qquad |V(z)|=\frac{||f'||_{\overline D}}{r^m}.$$

Evaluating $f'$ at $\alpha_i$ yields
$$\prod_{i=1}^{d}|f'(\alpha_i)|=\frac{(||f'||_{\overline D})^d}{r^{md}}\prod_{i=1}^{d}\prod_{j=1}^{m}|\alpha_i-c_j|, $$
and evaluating $f(z)-z_0$ at $c_j$ gives
$$\prod_{j=1}^{m}|f(c_j)-z_0|=\frac{R^m}{r^{dm}}\prod_{j=1}^{m}\prod_{i=1}^{d}|c_j-\alpha_i|,$$
where if $m=0$, products take value $1$.
It follows that 
$$\prod_{i=1}^{d}|f'(\alpha_i)|=(||f'||_{\overline D})^d\prod_{j=1}^{m}\frac{|f(c_j)-z_0|}{R}.$$
 The conclusion holds.
\end{proof}

The following result states a bound for the derivative on a pullback disk containing no critical point. 

\begin{lemma}\label{lem:no-value}
Let $z_0\in\mathbb{C}_p$, and let $\overline D_1, \overline D_2$ be two closed disks in $\mathbb{C}_p$ with $z_0\in \overline D_1$. Let $f\in\mathbb{C}_p(z)$ be a rational map with $f(\overline D_1)=\overline D_2$ such that the number of critical points of $f$ in $\overline D$, counted with multiplicity, is at most $\deg (f|_{\overline D})$ for any closed disk $\overline D\subseteq \overline D_1$. Suppose that  $\overline B_2\subsetneq \overline D_2$ is a closed subdisk containing  no critical value of $f|_{\overline D_1}$ and $f(z_0)\not\in \overline B_2$. 
Choose $\eta$ with $0<\eta\le1$ such that, for every critical value $v$ of $f|_{\overline D_1}$, 
$$\dist(v,\overline B_2)\ge\eta\cdot\dist(f(z_0),\overline B_2).$$
Then there is a closed subdisk $\overline B_1\subsetneq\overline D_1$ containing no critical point of $f$ such that $f(\overline B_1)=\overline B_2$  and, for any $z\in\overline B_1$
$$|f'(z)|=||f'||_{\overline B_1}\ge\eta\cdot|f'(z_0)|.$$ 
\end{lemma}
\begin{proof}
Let $r:=\dist(f(z_0), \overline B_2)>0$. Since $f(z_0)\not\in \overline B_2$, we have that 
$$\overline B_2\subseteq \overline D(f(z_0), r)\subseteq\overline D_2.$$
 Let $\overline D_0\subseteq\overline D_1$ be the disk containing $z_0$ such that $f(\overline D_0)=\overline D(f(z_0), r)$. Let $d$ denote the degree of $f$ on $\overline D_0$, and let $k$ denote the number of critical points of $f$ in $\overline D_0$, counted with multiplicity. Then  $k\le d$ by the assumption. 

Now fix $y\in \overline B_2$.   Since $|v-y|\ge \dist(v,\overline B_2)\ge \eta r$ for any critical value $v$ of  $f|_{\overline D_1}$, by Lemma \ref{lem:any}, we conclude that 
$$\prod_{\alpha\in \overline{D}_0: f(\alpha)=y}|f'(\alpha)|=(||f'||_{\overline D_0})^d\prod_{c\in \overline D_0: f'(c)=0}\frac{|f(c)-y|}{r}\ge (||f'||_{\overline D_0})^d\eta^k\ge (\eta\cdot||f'||_{\overline D_0})^d.$$
Thus, there is $\alpha_\ast\in\overline{D}_0$ with $f(\alpha_\ast)=y$ such that 
$$|f'(\alpha_\ast)|\ge \eta\cdot||f'||_{\overline D_0}\ge\eta\cdot|f'(z_0)|.$$
Let $\overline B_1\subseteq\overline D_0$ be the disk containing $\alpha_\ast$ such that $f(\overline B_1)= \overline B_2$. Then $\overline B_1$ contains no critical point, since $\overline B_2$ contains no critical value of $f|_{\overline D_1}$. It follows from Lemma \ref{lem:basic} that $|f'(z)|=|f'(\alpha_\ast)|$ for all $z\in \overline B_1$, and hence the conclusion holds. 
\end{proof}

\section{Return maps and ordered multigraphs}
In this section, let $f\in\mathbb{C}_p(z)$ be a rational map of degree at least $2$. In  \S\ref{sec:puzzles}, we construct puzzle pieces from a given finite union of pairwise disjoint closed disks. Then in \S\ref{sec:nests}, we use these puzzles to construct suitable neighborhoods of critical points (cf. principal nests in \cite[\S3.1]{Lyubich97}) and discuss the associated return maps under the assumption that $J(f)$ contains an infinite, compact, forward-invariant, minimal set $M$ that contains at least one  critical point of $f$. Finally in \S\ref{sec:graph}, we study the  ordered collection of multigraphs induced by the return maps.

\subsection{Puzzles}\label{sec:puzzles}
Let $U$ be a finite union of pairwise disjoint closed disks with radius in $|\mathbb{C}_p^\times|$  such that the image under $f$ of each disk is a closed disk of diameter less than $1$. For $N\ge 0$, set 
$$U_N=\bigcap_{j=0}^Nf^{-j}(U)\ \ \text{and}\ \ K_U=\bigcap_{N\ge 0}U_N.$$
Observing that  $U_{N+1}=U\cap f^{-1}(U_N)$, and inductively, applying \cite[Lemma 14.8]{Benedetto19}, we have that $U_N$ is a finite union of pairwise disjoint closed disks. 
Moreover, for any $x\in K_U$, let $B_N(x)$ be the maximal disk in $U_N$ containing $x$, and write $r_N(x):=\mathrm{diam}(B_N(x))$ for the diameter of $B_N(x)$. Observe that the disks $B_N(x)$ form a  nested sequence.

\begin{lemma}\label{lem:puzzle}
With this notation, the following statements hold:
\begin{enumerate}
\item For $N\ge 1$ and $y\in K_U$, if $B_N(y)\subsetneq B_0(y)$, then $f(B_N(y))=B_{N-1}(f(y))$.
\item If $D\subseteq K_U$ is a disk having positive diameter, then $D\subseteq F(f)$.
\item If $y\in K_U\cap J(f)$, then $r_N(y)\to0$ as $N\to\infty$. Moreover, if $X\subseteq  K_U\cap J(f)$ is compact, this convergence is uniform on $X$.
\end{enumerate}
\end{lemma}
\begin{proof}
Note that $B_N(y)$ is the maximal disk in $B_0(y)\cap f^{-1}(B_{N-1}(f(y)))$ containing $y$. It follows that $f(B_N(y))\subseteq f(B_0(y))\cap B_{N-1}(f(y))$. Since $B_N(y)\subsetneq B_0(y)$, we have that  $B_{N-1}(f(y))\subsetneq f(B_0(y))$. 
Then by  the maximality of  $B_N(y)$ and \cite[Lemma 14.8]{Benedetto19}, we conclude  that $f(B_N(y))=B_{N-1}(f(y))$. Hence statement (1) holds.

For statement (2), note that for each $j\ge 0$, the image $f^j(D)$ is a disk  contained in one component of $U$. Since each disk in $U$ has image, under $f$, with radius less than $1$, by Lemma \ref{lem:basic}, we conclude that, for any $x,y\in D$ and any $j\ge 0$, 
$$|f^{j+1}(x)-f^{j+1}(y)|\le \frac{1}{\diam(D)}|x-y|.$$
Thus $D\subseteq F(f)$. 

For statement (3), the convergence $r_N(y)\to0$ follows immediately from statement (2). We now show the  uniformity on $X$. Fix $\epsilon>0$. For each $x\in X$, since $r_N(x)\to0$, we choose a disk $B_N(x)$ with $r_N(x)<\epsilon$. Such disks form an open cover of $X$. Since $X$ is compact, we can find a  finite subcover from these disks, each of which has diameter  $<\epsilon$. Thus, there is a uniform $N$ such that $r_N(x)<\epsilon$ for all $x\in X$. 
\end{proof}

\subsection{Return maps}\label{sec:nests}
Denote by $\Crit(f)$ the set of critical points of $f$ in $\mathbb{P}^1(\mathbb{C}_p)$. Assume that $M\subseteq J(f)\cap K$ is an infinite, compact, forward-invariant, minimal set such that
$M\cap\Crit(f)\not=\emptyset$ and $M$ contains no poles. Write
$$M\cap\Crit(f)=\{c_1,\dots,c_s\}$$
where $s\ge 1$. 
Choose pairwise disjoint small closed disks $V_i^0$ containing $c_i$ such that $\diam(f(V_i^0))<1$,  and set 
$$V^0=\bigcup_{i=1}^s V_i^{0}.$$

For $x\in V^{0}$, we define its \emph{level $0$ return time $\tau_0$} by 
$$\tau_0(x):=\min\{k\ge 1: f^k(x)\in V^0\},$$ 
with $\tau_0(x)$ left undefined if $f^k(x)\not\in V^0$ for all $k\ge 1$. Let $S_0\subseteq V^0$ denote the set where $\tau_0$ is well-defined. 
Since $M$ is minimal, we have 
$$A_0:=M\cap  V^0\subseteq S_0.$$
 Define $\Dom_0$ to be the union of all maximal disks in $S_0$ that intersect $A_0$. The \emph{level $0$ return map} $R_0: \Dom_0\to V^0$ is defined by 
$$R_0(x):=f^{\tau_0(x)}(x).$$
We say a maximal disk $D$ of $\Dom_0$ is a \emph{level 0 source} if the return time $\tau_0(x)$ is independent of $x\in D$. 
For each $1\le i\le s$, let $V_i^1$ denote the source of $\Dom_0$ containing $c_i$.

We now proceed inductively. For each $n\ge 1$, define 
$$V^n=\bigcup_{i=1}^s V_i^n.$$
For $x\in V^n$, define its \emph{level $n$ return time}
$$\tau_n(x):=\min\{k\ge 1: f^k(x)\in V^n\},$$ 
with $\tau_n(x)$ left undefined if $f^k(x)\not\in V^n$ for all $k\ge 1$.
As before,  let $S_n\subseteq V^n$ denote the set on which $\tau_n$ is well-defined, and set $\Dom_n$ to be the union of all maximal disks in $S_n$ that intersect $A_n:=M\cap V^n\subseteq S_n$. Define the corresponding \emph{level $n$ return map} by 
$$R_n(x):=f^{\tau_n(x)}(x).$$
A \emph{level $n$ source} is a maximal disk in $\Dom_n$ on which the level $n$ return time is constant.
For each $1\le i\le s$, let $V_i^{n+1}$ denote the source in $\Dom_n$ containing $c_i$. 
Observe that,  for $x\in\Dom_{n+1}\subseteq\Dom_n$, 
$$\tau_{n+1}(x)\ge \tau_{n}(x).$$

 The restriction of the return map $R_n: V^{n+1}\to V^n$ is called \emph{central}, and the restriction $R_n: \Dom_n\setminus V^{n+1}\to V^n$ is called \emph{non-central}. 

\begin{lemma}\label{lem:nest}
There are $V^{0}$ and maps $\sigma_n:\{1,\dots,s\}\to\{1,\dots,s\}$ for $n\ge 0$ such that the following hold:
\begin{enumerate}
\item  $R_n(V_i^{n+1})=V_{\sigma_n(i)}^n$ and  $V_i^{n+1}\cap\Crit(R_n)=\{c_i\}$  for each $1\le i\le s$. 
\item If $D\subseteq\Dom_n$ is a level $n$ source with $D\cap V^{n+1}=\emptyset$, then $D\cap\Crit(R_n)=\emptyset$ and $R_n(D)=V_j^n$ for some $1\le j\le s$. 
\item For each $1\le i\le s$, $\diam(V_i^n)\to 0$, as $n\to\infty$. 
\end{enumerate}
\end{lemma}
\begin{proof}
Since $M$ is compact, choose a finite  union $U$ of disjoint closed disks with radius in $|\mathbb{C}_p^\times|$ covering $M$ such that $U\cap\Crit(f)\subseteq M$ and that each disk contains at most one critical point counted without multiplicity. Since $M$ is forward-invariant, we have that $M\subseteq K_U$. By Lemma \ref{lem:puzzle} (3), there exists an integer  $N\ge 0$ such that, for any $x\in M$ and  any $\ell\ge N$,  the disk $B_\ell(x)$ is a proper subdisk of the component of $U$ containing $x$. Let $V_i^0=B_N(c_i)$ for $1\le i\le s$. Then such disks $V_i^0$ are disjoint and $V_i^0\cap\Crit(f)=\{c_i\}$. 

Since $M$ is minimal, for each $c_i$, there exists a smallest integer $n_i\ge 1$ such that 
$$f^{n_i}(c_i)\in V^0:=\bigcup_{j=1}^sV_j^0.$$
 Define $V_i^1:=B_{N+n_i}(c_i)$. Then by Lemma \ref{lem:puzzle} (1), 
$$R_0(V_i^1)=f^{n_i}(V_i^1)=f^{n_i}(B_{N+n_i}(c_i))=B_N(f^{n_i}(c_i)).$$
Since $R_0(V_i^1)\subseteq V^0$, there is $1\le i_\ast\le s$ such that 
$$B_N(f^{n_i}(c_i))=B_N(c_{i_\ast})=V_{i_\ast}^0.$$
Then we define $\sigma_0(i)=i_\ast$. Note that $f^j(B_{N+n_i}(c_i))=B_{N+n_i-j}(f^j(c_i))$ by Lemma \ref{lem:puzzle} (1), which implies that,  for $1\le j\le n_i-1$, 
$$f^j(B_{N+n_i}(c_i))\cap\Crit(f)=\emptyset,$$ 
since $B_{N+n_i-j}(f^j(c_i))\cap V^0=\emptyset$. Noting that $B_{N+n_i}(c_i)\cap\Crit(f)=\{c_i\}$,  we conclude that  $V_i^1\cap\Crit(R_0)=\{c_i\}$. 

If $D\subseteq\Dom_0$ is a source with $D\cap V^1=\emptyset$, pick  $x\in D\cap A_0$ and choose $1\le j\le s$ such that $f^{\tau_0(x)}(x)\in V_j^0$. Then 
$$D=B_{N+\tau_0(x)}(x).$$ 
Indeed, by  Lemma \ref{lem:puzzle} (1), we have $B_{N+\tau_0(x)}(x)\subseteq D$, and, if $B_{N+\tau_0(x)}(x)\subsetneq D$, then there exists $y\in D\setminus B_{N+\tau_0(x)}(x)$  with $|y-x|/\diam(B_{N+\tau_0(x)}(x))$ sufficiently close to $1$ such that $\tau_0(y)>\tau_0(x)$ since $f^{\tau_0(x)}(y)\not\in V^0$, which is a contradiction.  Again by  Lemma \ref{lem:puzzle} (1), we have that  $D\cap\Crit(R_0)=\emptyset$ and $R_0(D)=B_N(f^{\tau_x}(x))= V_j^0$. 

Thus statements (1) and (2) hold for $n=0$. For $n\ge 1$, assume that statements (1) and (2) hold for $n-1$. Then for $1\le j\le s$, there is  $N_j>N$ such that $V_j^n=B_{N_j}(c_j)$.  For each $c_i$, there exists a smallest integer $n'_i\ge 1$ such that 
$$f^{n'_i}(c_i)\in V^n:=\cup_{j=1}^sV_j^n.$$
Define $V_i^{n+1}:=B_{N_i+n'_i}(c_i)$.
Then by Lemma \ref{lem:puzzle} (1), 
$$R_n(V_i^{n+1})=f^{n'_i}(V_i^{n+1})=f^{n'_i}(B_{N_i+n'_i}(c_i))=B_{N_i}(f^{n'_i}(c_i)).$$
Since $R_n(V_i^{n+1})\subseteq V^n$, there exists  $1\le i'_\ast\le s$ such that 
$$B_{N_i}(f^{n'_i}(c_i))=V_{i'_\ast}^n.$$ 
We define $\sigma_n(i)=i'_\ast$.
Note that $f^j(B_{N_i+n'_i}(c_i))=B_{N_i+n'_i-j}(f^j(c_i))$ by Lemma \ref{lem:puzzle} (1), which implies that,  for $1\le j\le n'_i-1$,  
$$f^j(B_{N_i+n'_i}(c_i))\cap\Crit(f)=\emptyset,$$ 
since $B_{N_i+n'_i-j}(f^j(c_i))\cap V^n=\emptyset$. Noting that $B_{N_i+n'_i}(c_i)\cap\Crit(f)=\{c_i\}$,  we conclude that  $V_i^{n+1}\cap\Crit(R_n)=\{c_i\}$. 

If $D\subseteq\Dom_{n}$ is a source with $D\cap V^{n+1}=\emptyset$, pick  $x\in D\cap A_n$ and let $1\le i, j\le s$ such that $x\in V_i^n$ and $f^{\tau_n(x)}(x)\in V_j^n$. Then 
$$D=B_{N_j+\tau_n(x)}(x).$$
 Indeed, by Lemma \ref{lem:puzzle} (1), we have $B_{N_j+\tau_n(x)}(x)\subseteq D$, and  moreover, there is $y\in D\setminus B_{N_j+\tau_n(x)}(x)$  with $|y-x|/\dim(B_{N_j+\tau_n(x)}(x))$ sufficiently closes to $1$ such that $\tau_n(y)>\tau_n(x)$ since $f^{\tau_n(x)}(y)\not\in V^n$. Again by  Lemma \ref{lem:puzzle} (1), we have that  $D\cap\Crit(R_n)=\emptyset$ and $R_n(D)=B_{N_j}(f^{\tau_n}(x))= V_j^n$. 

Thus we obtain statements (1) and (2) for $n$. Finally, note that $V_i^n\subseteq B_{N+n}(c_i)$. Then  $\diam(V_i^n)\to 0$ by Lemma \ref{lem:puzzle} (3). Thus statement (3) holds. 
\end{proof}

\subsection{Ordered multigraphs}\label{sec:graph}
In this subsection, we retain the notion in \S\ref{sec:nests} and consider the return maps $\{R_n\}_{n\ge 0}$ in Lemma \ref{lem:nest}. We construct a natural sequence of multigraphs from the return maps and analyze the nontrivial closed walks in these multigraphs. 
\begin{definition}
An \emph{ordered collection of multigraphs} induced by the return maps $\{R_n\}_{n\ge 0}$ is a collection of the directed multigraphs $\mathcal{G}_n$ with vertices $M\cap\Crit(f)=\{c_1,\dots,c_s\}$ such that the edges of $\mathcal{G}_n$  have two types: central edges and non-central edges, induced by the return maps; more precisely, there is a central edge from $c_i$ to $c_j$ of  $\mathcal{G}_n$ if and only if $R_n(V_i^{n+1})=V_j^n$, equivalently, $\sigma_n(i)=j$; and there is a non-central edge from $c_i$ to $c_j$ of  $\mathcal{G}_n$  if and only if there is a source $D\subseteq V_i^n\setminus V_i^{n+1}$ such that $R_n(D)=V_j^n$. The number of non-central edges from $c_i$ to $c_j$ is exactly the number of  distinct sources $D$. 
\end{definition}

Observe that $\mathcal{G}_n$ has finitely many vertices, and each vertex has exactly one outgoing central edge. Thus $\mathcal{G}_n$ has a nontrivial closed walk. 
 Let $\mathcal{N}_n$ be the non-central subgraph of $\mathcal{G}_n$ obtained by removing the central edges from $\mathcal{G}_n$. We establish the following two results,  depending on whether $\mathcal{N}_n$ contains a closed walk. 

\begin{proposition}\label{prop:main1}
Suppose there exists  $n\ge 0$ such that $\mathcal{N}_n$ has a nontrivial closed  walk. Then $\Rep(f)\not=\emptyset$.
\end{proposition}

\begin{proposition}\label{prop:main2}
If $\mathcal{N}_n$ has no nontrivial closed walk for any $n\ge 0$, then there exist integers $\mu\ge 1$, $n_0\ge 1$ and a list of distinct integers $\{i_j:0\le j\le\mu-1\}\subseteq\{1,\dots,s\}$ such that for any $n> n_0$, the graph $\mathcal{G}_n$ has a closed walk consisting of central edges from $c_{i_j}$ to $c_{i_{j+1}}$ for all $0\le j\le\mu-1$,  where the subscripts are modulo $\mu$. 
\end{proposition}

To prove Proposition \ref{prop:main1}, we first give some sufficient conditions for the existence of repelling periodic points.  

\begin{lemma}\label{lem:main}
Fix $n\ge 0$ and $1\le i, j\le s$. Suppose there exist a closed disk $\overline D\subseteq V_i^n$ and an integer $N\ge 1$ such that 
\begin{enumerate}
\item $\overline D\cap M=\emptyset$,
\item  $\overline D\cap \Crit(f^N)=\emptyset$,
\item  $f^N(\overline D)=V_j^n$. 
\end{enumerate}
Then  $\Rep(f)\not=\emptyset$. 
\end{lemma}

\begin{proof}
Since $\overline D\cap\Crit(f^N)=\emptyset$, we have that $|(f^N)'|$ is constant on $\overline D$ by Lemma \ref{lem:basic}. 
Set 
$$\Lambda:=||(f^N)'||_{\overline D},\ r=\dist(c_i,\overline D), \ \epsilon_0=\min_{z\in M\cap V_i^n}\dist(z,\overline D),\  \eta=\min\{1, \epsilon_0/r\}.$$
Since $M$ is compact and minimal, we can find an integer $m_0>0$ such that for any $z\in M$, there is  $0\le m_z\le m_0$ with $f^{m_z}(z)\in V_j^n$. 
Then choose 
$$C\ge \max\{|(f^q)'(z)|: z\in M, 0\le q\le m_0\}.$$
Consider the return map 
$$R_n|_{V_i^{n+1}}=f^{\tau_n(c_i)}|_{V_i^{n+1}}:V_i^{n+1}\to V_{\sigma_n(i)}^n.$$ 
Since $(R_n|_{V_i^{n+1}})'(c_i)=0$, we can choose a closed disk $\overline B\subseteq V_i^{n+1}$ containing $c_i$ such that 
$$\diam(\overline B)<r \ \ \text{and} \ \ ||(R_n|_{V_i^{n+1}})'||_{\overline B}\le \Lambda\eta/C.$$

Now fix $x\in M\cap \overline B$, and set $y=R_n|_{V_i^{n+1}}(x)$. Since $M$ is forward-invariant and minimal,  let $0\le m_y\le m_0$  be the smallest integer such that $u:=f^{m_y}(y)\in V_j^n$, and  let $\ell_u$ be the  smallest integer such that $f^{\ell_u}(u)\in  \overline B\subseteq V_i^n$. Consider the sequence of all components $\{V_{j_k}^n\}_{k=0}^{k_0}$  of $V^n$ containing a point in $\{f^\ell(u)\}_{\ell=0}^{\ell_u}$ with $V_{j_0}^n=V^n_j$ and $V_{j_{k_0}}^n=V_i^n$ such that that there is a level $n$ source contained in $V_{j_k}^n$ and containing some $f^\ell(u)$ returns to $V_{j_{k+1}}^n$ when this source returns to $V^n$. Successively pulling back $V_i^n$ along these returns, we obtain a closed disk $\Omega\subseteq V_j^n$ containing $u$ such that 
$$f^{\ell_u}(\Omega)=V_i^n.$$
By Lemma \ref{lem:nest} (1) and (2),  each  critical point of $f^{\ell_u}|_\Omega$ is an iterated preimage of one of the critical points $c_1,\dots,c_s$.  
Counting the critical points with multiplicity and the degrees, we have that, for any closed disk $\Omega'\subseteq\Omega$, 
$$|\Crit(f^{\ell_u})\cap \Omega'|\le\deg(f^{\ell_u}|_{\Omega'})-1,$$
Moreover, by the  forward-invariance of $M$, every critical value of $f^{\ell_u}|_\Omega$ is contained in $M\cap V_i^n$. 

Note that $\overline D$ contains no critical value of $f^{\ell_u}|_\Omega$ since $\overline D\cap M=\emptyset$. Moreover, since  $f^{\ell_u}(u)\in  \overline B$ and $\diam(\overline B)<r$, we have that 
$$\dist(f^{\ell_u}(u), \overline D)=\dist(c_i,\overline D)=r,$$
 and hence for every critical value $v$ of $f^{\ell_u}|_\Omega$, 
$$\dist(v,\overline D)\ge\epsilon_0\ge\eta r=\eta\cdot\dist(f^{\ell_u}(u),\overline D).$$ 
Thus by Lemma \ref{lem:no-value},  there exists a closed subdisk $\Omega_\ast\subsetneq\Omega$ containing no critical point of $f^{\ell_u}$ such that $f^{\ell_u}(\Omega_\ast)=\overline D$  and, for any $z\in\Omega_\ast$, 
$$\Theta:=||(f^{\ell_u})'||_{\Omega_\ast}=|(f^{\ell_u})'(z)|\ge\eta\cdot|(f^{\ell_u})'(u)|.$$ 

Now suppose, for contradiction, that $\Rep(f)=\emptyset$. Consider the restriction $f^{N+\ell_u}:\Omega_\ast\to V_j^n$, which is surjective and has no critical point. Note that, for any $z\in \Omega_\ast$, 
$$|(f^{N+\ell_u})'(z)|=\Lambda\Theta.$$   
 If $\Omega_\ast\subsetneq V_j^n$, by Lemma \ref{lem:fix}, the absence of repelling periodic points of $f$ implies that $\Lambda\Theta\le1$. If  $\Omega_\ast=V_j^n$, then we also have $\Lambda\Theta\le1$ by Lemma \ref{lem:basic}. Then setting 
 $$\kappa(x):=\tau_n(c_i)+m_y+\ell_u,$$ 
 we have 
$$|(f^{\kappa(x)})'(x)|=|(f^{\ell_u})'(u)|\cdot|(f^{m_y})'(y)|\cdot|(R_n|_{V_i^{n+1}})'(x)|\le|(f^{\ell_u})'(u)|\cdot C \cdot (\Lambda\eta/C)\le \Lambda\Theta\le1.$$ 

Since again $M$ is compact and minimal,  we can find an integer $\ell_0>0$ such that for any $z\in M$, there is $0\le \ell_z\le\ell_0$ with  $f^{\ell_z}(z)\in\overline B$. Then consider the small closed disks in $\overline B$ such that on each disk, the corresponding $m_y$ and $\ell_u$ are constant. We can choose finitely many such disks $W_1,\dots, W_{k_0}$ to cover $M\cap\overline B$. For each $1\le j\le k_0$, set $\kappa_j=\kappa(w)$ for any $w\in M\cap W_j$ and denote $F_j:=f^{\kappa_j}$. Then $F_j(W_j)\subseteq\overline B$ and,  for any $w\in M\cap W_j$, 
$$|F_j'(w)|=|(f^{\kappa_j})'(w)|\le \Lambda\Theta\le1.$$ 
Moreover, let 
$$\kappa_{\max}:=\max\{\kappa_j:1\le j\le k_0\},$$
and if necessary, divide each $W_j$ into smaller closed disks so that $f^q$ has no poles in $W_j$ for any $1\le q\le \kappa_{\max}$, since $M$ is forward-invariant and contained in $K$. 

Let $\delta:=\diam(\overline B)$ and let $\alpha:=\min\{\diam(W_j):1\le j\le k_0\}$. Choose small $0<\epsilon\le\min\{\alpha,\delta\}$ such that
$$\delta\epsilon/\alpha^2\le 1.$$
For $1\le j\le k_0$,  choose $x_0\in M\cap W_j$, and write 
$$F_j(x_0+z)=F_j(x_0)+\sum_{q\ge 1}a_{j,q}(x_0)z^q.$$
Since $\overline D(x_0,\alpha)\subseteq W_j$, we have that 
$$F_j(\overline D(x_0,\alpha))\subseteq F_j(W_j)\subseteq\overline B.$$
 It follows that, for any $q\ge 1$, 
 $$|a_{j,q}(x_0)|\alpha^q\le \delta.$$
  Thus, for $q\ge 2$, 
$$|a_{j,q}(x_0)|\epsilon^q\le \frac{\delta\epsilon^q}{\alpha^q}= \frac{\delta\epsilon}{\alpha^2}\left(\frac{\epsilon}{\alpha}\right)^{q-2}\epsilon\le\epsilon.$$
For $q=1$,  since $|F_j'(x_0)|\le1$, we have 
$$|a_{j,1}(x_0)|\epsilon\le \epsilon.$$
Thus for any $w\in M\cap\overline B$,  letting $1\le j\le k_0$ be such that $w\in W_j$, we have  $F_j(w)\in M\cap\overline B$ and $F_j$ maps $\overline{D}(w,\epsilon)$ into $\overline{D}(F_j(w),\epsilon)\subseteq\overline B$. Thus, $F_j$ is $1$-Lipschitz on $\overline{D}(w,\epsilon)$ by Lemma \ref{lem:basic}.  

 Again by Lemma \ref{lem:basic}, let $L>0$ be a number greater than all Lipschitz constants for $f^q$ on $W_j$ for $1\le q\le \kappa_{\max}$ and $1\le j\le k_0$. We conclude that $f^q$ is $L$-Lipschitz on $\overline{D}(c_i,\epsilon)$ for any $q\ge 1$. Indeed, we consider integers $q_0=0, q_1,\cdots, q_m\ge 1$ defined by  
 $$q_h:=\kappa(f^{\sum_{k=0}^{h-1}q_k}(c_i))$$ 
 for $1\le h\le m$ and 
 $$0\le q-\sum_{k=1}^mq_k<\kappa(f^{\sum_{k=1}^{m}q_k}(c_i)),$$  and then noting that  $f^{\sum_{k=1}^{m}q_k}$ is $1$-Lipschitz on $\overline{D}(c_i,\epsilon)$  and that $0\le q-\sum_{k=1}^mq_k\le \kappa_{\max}$, we have that  $f^q$ is $L$-Lipschitz on $\overline{D}(c_i,\epsilon)$. This  contradicts $c_i\in J(f)$. 
\end{proof}

\begin{proof}[Proof of Proposition \ref{prop:main1}]

Since $\mathcal {N}_n$ has a nontrivial closed walk, there exist $1\le i\le s$, $m\ge 1$ and a closed disk $\overline D\subseteq V_i^n$ such that $f^m(\overline D)=V_i^n$, and that $\overline D\cap\Crit(f^m)=\emptyset$ by Lemma \ref{lem:nest} (2) since each directed edge in this closed walk corresponds to a non-central return on a source. For $j\ge 1$, set 
$$W_j=\{z\in\overline D: f^{km}(z)\in\overline D\ \ \text{for all}\  \ 0\le k<j\}.$$

Note that $W_j$ is a nonempty finite union of closed disks, each component of which maps onto $V_i^n$ under $f^{jm}$ without critical points. Moreover, $W_{j+1}\subseteq W_j$.
 If every $W_j$ intersects $M$, compactness gives a point 
 $$x\in M\cap(\bigcap_j W_j).$$ 
 Since $M$ is minimal, the forward orbit of $x$ is dense in $M$. However, at returns,  this orbit follows non-central source, ad hence it is disjoint from $V^{n+1}$, which is a contradiction. 
 
Thus there is $j_0\ge 1$ such that  $W_{j_0}\cap M=\emptyset$. Let $D_0$ be a maximal disk in $W_{j_0}$. Then $f^{j_0m}(D_0)=V_i^n$. Applying Lemma \ref{lem:main}  to $D_0$ and $f^{j_0m}$, we conclude that $\Rep(f)\not=\emptyset$. 
\end{proof}

To prove Proposition \ref{prop:main2}, we first establish the following lemma which upgrades certain paths to non-central edges by increasing the levels of the graphs. 
\begin{lemma}\label{lem:connect}
For $1\le i, j\le s$, if $\mathcal{N}_n$ has a path from $c_{\sigma_n(i)}$ to $c_j$ with $j\not=\sigma_{n+1}(i)$, then $\mathcal{N}_{n+1}$ has an edge from $c_i$ to $c_j$. 
\end{lemma}
\begin{proof}
Pull back $V_j^{n+1}$ along a path $\mathcal{C}$  from $c_{\sigma_n(i)}$ to $c_j$ in $\mathcal{N}_n$, and then pull back by the map $R_n|_{V_i^{n+1}}$. We obtain a closed disk $\overline{D}\subseteq V_i^{n+1}$ such that $f^m(\overline{D})=V_j^{n+1}$ for some $m\ge 1$. Since $\mathcal{C}$ consists of noncentral edges, we have that $f^h(\overline D)\cap V^{n+1}=\emptyset$ for any $1\le h\le m-1$. Hence $f^m|_{\overline D}=R_{n+1}|_{\overline D}$.  Moreover, $c_i\notin\overline{D}$ since $j\not=\sigma_{n+1}(i)$. Thus $\mathcal{G}_{n+1}$ has a noncentral edge from $c_i$ to $c_j$. 
\end{proof}

\begin{proof}[Proof of Proposition \ref{prop:main2}]
For a closed walk $\mathcal{C}$ in $\mathcal{G}_n$, let $\nu_n(\mathcal{C})$ count its central edge
occurrences. Set
$$ \mu:=\min\{\nu_n(\mathcal{C}): n\ge0,\ \mathcal{C}\ \text{is a closed walk in}\ \mathcal{G}_n\}.$$
Since $\mathcal{N}_n$ has no closed  walk for all $n\ge 0$, we have that  $\mu\ge 1$.   Let $\mathcal{C}_\ast$ be a closed walk in $\mathcal{G}_{n_0}$ for some $n_0\ge 0$ such that $\nu_{n_0}(\mathcal{C}_\ast)=\mu$.

Cut $\mathcal{C}_\ast$ exactly before the central occurrences. Then each piece in this cut contains exactly one central edge in $\mathcal{G}_{n_0}$. Pick such a piece, say  from $c_i$ to $c_j$. If $j=\sigma_{n_0+1}(i)$, replace this piece by the central edge from $c_i$ to $c_j$ in $\mathcal{G}_{n_0+1}$; and  if  $j\not=\sigma_{n_0+1}(i)$,  by  Lemma \ref{lem:connect}, we  replace this piece by a non-central edge from $c_i$ to $c_j$ in $\mathcal{G}_{n_0+1}$. Then we obtain a closed walk $\mathcal{C}'_\ast$ which has exactly $\mu$ edges, each from $c_{i_j}$ to $c_{i_{j+1}}$ with $0\le j\le\mu-1$ and $c_{i_\mu}=c_{i_0}$.  Note that $\nu_{n_0+1}(\mathcal{C}'_\ast)\le \nu_{n_0}(\mathcal{C}_\ast)$. By the choice of $\mathcal{C}_\ast$, we conclude that 
$$\nu_{n_0+1}(\mathcal{C}'_\ast)=\nu_{n_0}(\mathcal{C}_\ast)=\mu.$$ 
Thus each piece in the cut of $\mathcal{C}_\ast$ is replaced by a central edge in  $\mathcal{G}_{n_0+1}$. It follows that the graph $\mathcal{G}_{n_0+1}$ has central edges from $c_{i_j}$ to $c_{i_{j+1}}$ for any $0\le j\le\mu-1$. 

Then working with $\mathcal{C}'_\ast$ and applying the above argument inductively, we conclude that for any $n>n_0$,  the graph $\mathcal{G}_n$ has central edges from $c_{i_j}$ to $c_{i_{j+1}}$ for any $0\le j\le\mu-1$. 
\end{proof}

\section{Proof of Theorem \ref{mainB}}
In this section, we prove Theorem \ref{mainB}. We retain the notation of the preceding section. 
\begin{proof}[Proof of Theorem \ref{mainB}]
Fix a non-empty compact forward-invariant set $X\subset J(f)$. Every inclusion chain of nonempty forward-invariant  compact subsets of $X$ has a non-empty forward-invariant intersection. Zorn's Lemma gives a nonempty, compact, forward-invariant and minimal set $M$ in $J(f)$. If $M$ is finite, then $M$ is a repelling periodic cycle, which finishes the proof. Now we assume that $M$ is infinite. 

\textbf{Case 1:} $M\cap \Crit(f)=\emptyset$. For each $x\in M$, consider a disk $D(x)$ such that $f$ has Weierstrass degree $1$ on $D(x)$ and $\diam(f(D(x)))<1$. Then such disks cover $M$. Since $M$ is compact, we can get a finite finite union $U$  of disjoint closed disks covering $M$.

Since $M$ is forward-invariant, we have that  $M\subseteq K_U$.  Then by Lemma \ref{lem:puzzle} (3), 
there is $N>0$ such that, for every $j\ge N$ and $y\in M$, the disk $B_j(y)$ is proper in $B_0(y)$. Pick $x\in M$. Since $M$ is minimal, we have that $x$ is recurrent. Hence there is $n:=n_x\ge1$ such that
$f^n(x)\in B_N(x)$. Then by  Lemma \ref{lem:puzzle} (1), we conclude that for $0\le i\le n$, 
$$f^i(B_{N+n}(x))=B_{N+n-i}(f^i(x))\subsetneq B_0(f^i(x)).$$
In particular, 
$$f^n(B_{N+n}(x))=B_{N}(f^n(x))=B_{N}(x).$$
Thus $f^n:B_{N+n}(x)\to B_{N}(x)$ has Weierstrass degree $1$. 

Since $x\in M\subseteq K_U\cap J(f)$,  Lemma \ref{lem:puzzle} (3) implies that $B_{N+n}(x)\subsetneq B_{N}(x)$. Thus $f^n$ has a fixed point $x_0\in B_{N+n}(x)$ with 
$$|(f^n)'(x_0)|=\frac{\diam(B_{N}(x))}{\diam(B_{N+n}(x))}>1.$$
Thus $x_0\in\Rep(f)$. 

\textbf{Case 2:} $M\cap \Crit(f)\not=\emptyset$. Write $M\cap \Crit(f)=\{c_1,\dots,c_s\}$ with $s\ge 1$. Consider the ordered collection of graphs $\{\mathcal{G}_n\}$ constructed in \S\ref{sec:graph}. 
If each nontrivial closed walk of $\mathcal{G}_n$ has a central edge for every $n$, then  by Proposition \ref{prop:main2}, there exist  distinct indices $1\le i_0,\ldots,i_{\mu-1}\le s$ and an integer $n_0$ such
that $\sigma_n(i_j)=i_{j+1}$ with $c_{i_\mu}=c_{i_0}$ for any $n>n_0$. It follows that $R_n(c_{i_j})\in V_{i_{j+1}}^n$  for $n>n_0$. Since the return times $\tau_n(c_{i_j})\le \tau_{n+1}(c_{i_j})$, we claim that 
$$\tau_n(c_{i_j})=\tau_{n+1}(c_{i_j}).$$
 For otherwise, since $f^{\tau_n(c_{i_j})}(c_{i_j})\in  V_{i_{j+1}}^n$, we have $f^{\tau_n(c_{i_j})}(c_{i_j})\notin  V_{i_{j+1}}^{n+1}$, and then obtain a non-central return $f^{\tau_{n+1}(c_{i_j})-\tau_n(c_{i_j})}$ from $D:=f^{\tau_n(c_{i_j})}(V_{i_j}^{n+2})\subseteq V_{i_{j+1}}^n\setminus V_{i_{j+1}}^{n+1}$ to $ V_{i_{j+1}}^{n+1}$, 
 which  induces a return from a disk $D'\subseteq V_{i_{j+1}}^n\setminus V_{i_{j+1}}^{n+1}$ containing $D$ to $V_{i_{j+1}}^n$ and hence gives a non-central closed walk in $\mathcal{G}_N$, contradicting the assumption that each non-trivial closed walk of $\mathcal{G}_n$ has a central edge. 
Since $\diam(V_i^n)\to 0$ as $n\to\infty$ for all $1\le i\le s$ by Lemma \ref{lem:nest} (3), we conclude that 
$$f^{\sum_{j=0}^{\mu-1}\tau_n(c_{i_j})}(c_{i_0})=c_{i_0}.$$
This contradicts $c_{i_0}\in J(f)$. Thus $\mathcal{G}_n$ has a closed walk containing no central edge for some $n$. Then by Proposition \ref{prop:main1}, $\Rep(f)\not=\emptyset$. 
 \end{proof}

 \section{Proof of Theorem \ref{main}}
 In this section, we deduce Theorem  \ref{main} from Theorem  \ref{mainB} and Proposition \ref{prop:algebraic} below.
 
 \begin{proposition}\label{prop:algebraic}
Let $f\in\overline{\mathbb{Q}}_p(z)$ be a rational map of degree at least two with $J(f)\not=\emptyset$. Then $J(f)\cap\overline{\mathbb{Q}}_p\not=\emptyset$.
 \end{proposition}
 
 The proof of Proposition \ref{prop:algebraic} combines the following two results via an argument by contradiction. The first result (Lemma \ref{lem:U}) provides an arithmetic disk partition, whereas the second result (Lemma \ref{lem:finite}) produces a repelling periodic point. 
  \begin{lemma}\label{lem:U}
Let $K/\mathbb{Q}_p$ be a finite extension, and let $f\in K(z)$ be a rational map of degree at least $2$ such that $J(f)$ is nonempty and contained in $D^-(0,1)$. Assume that all poles and critical points of $f$ are  contained in $\mathbb{P}^1(K)$.  If $\mathbb{P}^1(K)\subseteq F(f)$, then the collection $\mathcal{U}$ of all maximal open disks, intersecting $J(f)$,  in $D^-(0,1)\setminus K$ satisfies the following: 
 \begin{enumerate}
 \item there is $\epsilon>0$ such that $\epsilon\le\diam(D)\le 1$ for any $D\in\mathcal{U}$,
 \item if $D\in\mathcal{U}$ with $f(D)\cap K=\emptyset$, then $f(D)\in\mathcal{U}$,
 \item there are only finitely many disks $D\in\mathcal{U}$ such that $f(D)$ intersects $K$;  
 \end{enumerate}
 \end{lemma}
 \begin{proof}
Since $K\cap\overline D(0,1)$ is compact and contained in the open set $F(f)$, choose $\epsilon>0$ such that the $\epsilon$-neighborhood of $K\cap\overline D(0,1)$ is also contained in $F(f)$. Then  for any $x\in J(f)$ and any $a\in K\cap\overline D(0,1)$, 
 $$\epsilon\le |x-a|\le 1.$$
Thus statement $(1)$ holds.

For statement (2),  pick any $D\in\mathcal{U}$ with $f(D)\cap K=\emptyset$. By Lemma \ref{lem:basic}, the image $f(D)$ is a disk. The forward-invariance of $J(f)$ implies that $f(D)\cap J(f)\not=\emptyset$, and hence $f(D)\subsetneq D^-(0,1)\setminus K$. We claim that $f(D)$ is a maximal disk in $D^-(0,1)\setminus K$. Suppose otherwise. There exists a disk $B\subsetneq D^-(0,1)\setminus K$ such that $f(D)\subsetneq B$. Since $f(K)\subseteq \mathbb{P}^1(K)$, we have $f^{-1}(B)\cap K=\emptyset$. Let $B'$ be the maximal disk of $f^{-1}(B)$ containing $D$. Then $B'\subsetneq D^-(0,1)\setminus K$. Moreover,  by \cite[Proposition 3.29]{Benedetto19}, we have that $f(B')=B$. It follows that $D\subsetneq B'$, which contradicts the  maximality of $D$.

For statement (3), let $\mathcal{E}\subseteq\mathcal{U}$ be the set of  $D\in\mathcal{U}$ with $f(D)\cap K\not=\emptyset$. For any $D\in\mathcal{E}$, choose $y_D\in D$ such that $f(y_D)\in K$. Then 
$$y_D\in f^{-1}\bigl(\mathbb P^1(K)\bigr)\cap\overline D(0,1).$$
Moreover, for two distinct $D_1,D_2\in\mathcal{E}$, by statement (1), we have $|y_{D_1}-y_{D_2}|\ge\epsilon$. 
We claim that the set  $X:=f^{-1}\bigl(\mathbb{P}^1(K)\bigr)\cap\overline D(0,1)$ is compact. Indeed, for any $x\in f^{-1}\bigl(\mathbb P^1(K)\bigr)$, either $x$ is a pole or there is $y\in K$ such that $f(x)=y$. In either case, $[K(x):K]\le\deg f$. Since $K$ has only finitely many field extensions of degree at most $\deg f$ (see \cite[Chapter II \S5 Proposition 14]{Lang94}), let $L$ be the smallest subfield of $\mathbb{C}_p$ containing  these finite extensions. Then $L/K$ is a finite extension and $f^{-1}\bigl(\mathbb{P}^1(K)\bigr)\subseteq\mathbb{P}^1(L)$. It follows that $X$ is contained in the ring of integers  $\mathcal{O}_L$ of $L$. Note that $X$  is closed since $f$ is continuous, and that $\mathcal{O}_L$  is compact. We obtain that $X$ is compact. 
It follows that $X$ is a compact metric space, and hence totally bounded. We conclude that $\mathcal{E}$ is a finite set. 
\end{proof}

\begin{lemma}\label{lem:finite}
Let $f\in\mathbb{C}_p(z)$ be a rational map with nonempty $J(f)\subseteq D^-(0,1)$ such that $J(f)\cap\Crit(f)=\emptyset$. 
Suppose a collection $\mathcal{U}$ of pairwise disjoint open disks covers $J(f)$ and satisfies the following:
\begin{enumerate}
\item each disk in $\mathcal{U}$ contains neither  poles nor critical points of $f$;
\item there exists $\epsilon>0$ such that $\epsilon\le\diam(D)\le1$ for any $D\in \mathcal{U}$; 
\item for any $D_1,D_2$ in $\mathcal{U}$,  if $f(D_1)\cap D_2\not=\emptyset$, then $D_2\subseteq f(D_1)$;
\item $f(D)\in\mathcal{U}$  for all except finitely many  $D\in\mathcal{U}$. 
\end{enumerate}
Then $\Rep(f)\not=\emptyset$. 
\end{lemma}

\begin{proof}
Write $\lambda_\alpha=|f'(z)|$ for any $z\in D_\alpha$ and $D_\alpha\in\mathcal{U}$ by Lemma \ref{lem:basic}. Let $\mathcal{E}\subset\mathcal{U}$ be the finite set of disks in $\mathcal{U}$ whose image is not in $\mathcal{U}$, and denote $N:=|\mathcal{E}|$. 

Suppose, for contradiction, that $\Rep(f)=\emptyset$. Consider the directed graph $G$ whose vertices are the disks in $\mathcal{U}$ with an edge $D_\alpha\to D_\beta$ whenever $f(D_\alpha)\cap D_\beta\not=\emptyset$ for $D_\alpha, D_\beta\in\mathcal{U}$. Then  the absence of repelling periodic cycles implies that for every directed cycle $D_{\alpha_0}\to D_{\alpha_1}\to\cdots\to D_{\alpha_{n-1}}\to D_{\alpha_0}$ in $G$, 
$$\prod_{j=0}^{n-1}\lambda_{\alpha_j}\le1.$$

Assign to each edge $D_\alpha\to D_\beta$ in $G$ the weight
$$w_{\alpha\beta}:=\lambda_\alpha\frac{\diam(D_\alpha)}{\diam(D_\beta)},$$
and define 
$$\Theta:=\max\left(\{1\}\cup\left\{\frac{\lambda_\alpha\diam(D_\alpha)}{\epsilon}: D_\alpha\in\mathcal{E}\right\}\right).$$
Then for each edge $D_\alpha\to D_\beta$,  if $D_\alpha\not\in\mathcal{E}$, then $f(D_\alpha)=D_\beta$, and hence by Lemma \ref{lem:basic}, 
$$w_{\alpha\beta}\le 1.$$
If $D_\alpha\in\mathcal{E}$, then
$$w_{\alpha\beta}\le\Theta.$$

For a walk $\gamma: D_{\beta_0}\to D_{\beta_1}\to\cdots\to D_{\beta_{m-1}}\to D_{\beta_m}$ in $G$, define the weight of $\gamma$ by 
$$w(\gamma)=\prod_{j=0}^{m-1}w_{\beta_j\beta_{j+1}}.$$
 We claim that, for any walk $\gamma$ in $G$, 
$$w(\gamma)\le \Theta^N.$$
Indeed, noting that  for every directed cycle $D_{\alpha_0}\to D_{\alpha_1}\to\cdots\to D_{\alpha_{n-1}}\to D_{\alpha_0}$, 
$$\prod_{j=0}^{n-1}w_{\alpha_j\alpha_{j+1}}=\prod_{j=0}^{n-1}\lambda_{\alpha_j}\le1,$$
 where the subscripts of $\alpha_j$ are modulo $n$, 
 we remove all cycles in $\gamma$ and get a resulting path $\gamma'$ which visits each vertex in $\mathcal{E}$ at most once and has weight 
 $$w(\gamma)\le w(\gamma')\le\Theta^N.$$

Fix any $x\in J(f)$, and for any $j\ge 0$, let $D_{\alpha_j}$ be the unique disk in $\mathcal{U}$ containing $f^j(x)$. 
Since $\epsilon\le\diam(D_{\alpha_j})\le1$, we obtain that for any $n\ge 1$, 
$$|(f^n)'(x)|=\prod_{j=0}^{n-1}\lambda_{\alpha_j}=\frac{\diam(D_{\alpha_n})}{\diam(D_{\alpha_0})}\prod_{j=0}^{n-1}w_{\alpha_j\alpha_{j+1}}\le \frac{\Theta^N}{\epsilon}=:C.$$

To produce a contradiction, we show that the iterates $f^n$ are Lipschitz on a common neighborhood of $x$, with a constant independent of $n$. Set $r=\epsilon|p|^2$. 
For any $D_\alpha\in\mathcal{U}$ and  any $y\in D_\alpha$, write 
$$f(y+z)=f(y)+\sum_{\ell\ge1}b_\ell z^\ell.$$
Since $|f'|\equiv\lambda_\alpha$ on $D_\alpha$, we have that 
$$|\ell b_\ell|\left(\diam(D_\alpha)\right)^{\ell-1}\le\lambda_\alpha.$$
Then for \(\ell\ge2\) and  any $z\in\overline D(0,r)$, noting that  $|\ell|\ge |p|^{\ell-1}$, we conclude that, for $z\not=0$, 
$$|b_\ell||z|^\ell\le\frac{\lambda_\alpha r^{\ell-1}|z|}{|\ell|\left(\diam(D_\alpha)\right)^{\ell-1}}\le\frac{ \lambda_\alpha\epsilon^{\ell-1}|p|^{2(\ell-1)}|z|}{|\ell|\left(\diam(D_\alpha)\right)^{\ell-1}}\le\lambda_\alpha |p|^{\ell-1}|z|<\lambda_\alpha |z|.$$
Thus 
$$|f(y+z)-f(y)|=|b_1||z|=\lambda_\alpha|z|$$
Applying this identity inductively along  the  orbit of $x$ gives that, for any $n\ge 1$ and for any  $x_1,x_2\in\overline D(x,r/C)$,
$$|f^n(x_1)-f^n(x_2)|=|(f^n)'(x)|\,|x_1-x_2|\le C|x_1-x_2|.$$
Therefore, the iterates $f^n$ are uniformly Lipschitz on a neighborhood of $x$, which contradicts $x\in J(f)$.
\end{proof}

We now prove Proposition \ref{prop:algebraic}.
\begin{proof}[Proof of  Proposition \ref{prop:algebraic}]
Suppose, for contradiction, that $J(f)\cap\overline{\mathbb{Q}}_p=\emptyset$. 
We choose a finite extension $K$ of $\mathbb{Q}_p$  containing all coefficients, critical points and poles of $f$.  Conjugating $f$ by an element in $\mathrm{PGL}_2(K)$, we can assume that $J(f)\subseteq D^-(0,1)$. Since this conjugacy preserves algebraicity, we have $\mathbb{P}^1(K)\subseteq F(f)$.  Let $\mathcal{U}$ be the collection of maximal open disks, intersecting $J(f)$,  in $D^-(0,1)\setminus K$. Lemma \ref{lem:U} then shows that $\mathcal{U}$ satisfies the hypotheses  of Lemma \ref{lem:finite}, which implies that $J(f)$ contains a repelling periodic point. This contradicts the assumption that $J(f)\cap\overline{\mathbb{Q}}_p=\emptyset$. 
\end{proof}

We now prove Theorem \ref{main}.
\begin{proof}[Proof of Theorem \ref{main}]
If $J(f)=\emptyset$, then the conclusion holds immediately since $\Rep(f)\subseteq J(f)$. We now assume that $J(f)\not=\emptyset$. By Proposition \ref{prop:algebraic}, pick $x\in J(f)\cap \overline{\mathbb{Q}}_p$. Consider a finite extension $L$ of $\mathbb{Q}_p$ containing $x$ and all coefficients of $f$. Then the closure $M$ of the orbit of $x$ under $f$ is contained in $\mathbb{P}^1(L)$, and hence is nonempty, compact and forward-invariant. Noting that $M\subseteq J(f)$, we obtain $\Rep(f)\not=\emptyset$ by Theorem \ref{mainB}. 
\end{proof}

\bibliographystyle{siam}
\bibliography{references}	
\end{document}